\documentclass{amsart} 
\usepackage{latexsym} 
\usepackage{amssymb} 
\usepackage{amsmath} 

\begin{document}


\newtheorem{theorem}{Theorem} 
\newtheorem{problem}{Problem} 
\newtheorem{definition}{Definition} 
\newtheorem{lemma}{Lemma} 
\newtheorem{proposition}{Proposition} 
\newtheorem{corollary}{Corollary} 
\newtheorem{example}{Example} 
\newtheorem{conjecture}{Conjecture} 
\newtheorem{algorithm}{Algorithm} 
\newtheorem{exercise}{Exercise} 
\newtheorem{remarkk}{Remark} 
 
\newcommand{\be}{\begin{equation}} 
\newcommand{\ee}{\end{equation}} 
\newcommand{\bea}{\begin{eqnarray}} 
\newcommand{\eea}{\end{eqnarray}} 
\newcommand{\beq}[1]{\begin{equation}\label{#1}} 
\newcommand{\eeq}{\end{equation}} 
\newcommand{\beqn}[1]{\begin{eqnarray}\label{#1}} 
\newcommand{\eeqn}{\end{eqnarray}} 
\newcommand{\beaa}{\begin{eqnarray*}} 
\newcommand{\eeaa}{\end{eqnarray*}} 
\newcommand{\req}[1]{(\ref{#1})} 
 
\newcommand{\lip}{\langle} 
\newcommand{\rip}{\rangle} 

\newcommand{\uu}{\underline} 
\newcommand{\oo}{\overline} 
\newcommand{\La}{\Lambda} 
\newcommand{\la}{\lambda} 
\newcommand{\eps}{\varepsilon} 
\newcommand{\om}{\omega} 
\newcommand{\Om}{\Omega} 
\newcommand{\ga}{\gamma} 
\newcommand{\rrr}{{\Bigr)}} 
\newcommand{\qqq}{{\Bigl\|}} 
 
\newcommand{\dint}{\displaystyle\int} 
\newcommand{\dsum}{\displaystyle\sum} 
\newcommand{\dfr}{\displaystyle\frac} 
\newcommand{\bige}{\mbox{\Large\it e}} 
\newcommand{\integers}{{\Bbb Z}} 
\newcommand{\rationals}{{\Bbb Q}} 
\newcommand{\reals}{{\rm I\!R}} 
\newcommand{\realsd}{\reals^d} 
\newcommand{\realsn}{\reals^n} 
\newcommand{\NN}{{\rm I\!N}} 
\newcommand{\DD}{{\rm I\!D}} 
\newcommand{\degree}{{\scriptscriptstyle \circ }} 
\newcommand{\dfn}{\stackrel{\triangle}{=}} 
\def\complex{\mathop{\raise .45ex\hbox{${\bf\scriptstyle{|}}$} 
     \kern -0.40em {\rm \textstyle{C}}}\nolimits} 
\def\hilbert{\mathop{\raise .21ex\hbox{$\bigcirc$}}\kern -1.005em {\rm\textstyle{H}}} 
\newcommand{\RAISE}{{\:\raisebox{.6ex}{$\scriptstyle{>}$}\raisebox{-.3ex} 
           {$\scriptstyle{\!\!\!\!\!<}\:$}}} 
 
\newcommand{\hh}{{\:\raisebox{1.8ex}{$\scriptstyle{\degree}$}\raisebox{.0ex} 
           {$\textstyle{\!\!\!\! H}$}}} 

\newcommand{\OO}{\won} 
\newcommand{\calA}{{\mathcal A}} 
\newcommand{\calB}{{\cal B}} 
\newcommand{\calC}{{\cal C}} 
\newcommand{\calD}{{\cal D}} 
\newcommand{\calE}{{\cal E}} 
\newcommand{\calF}{{\mathcal F}} 
\newcommand{\calG}{{\cal G}} 
\newcommand{\calH}{{\cal H}} 
\newcommand{\calK}{{\cal K}} 
\newcommand{\calL}{{\mathcal L}} 
\newcommand{\calM}{{\cal M}} 
\newcommand{\calO}{{\cal O}} 
\newcommand{\calP}{{\cal P}} 
\newcommand{\calU}{{\mathcal U}} 
\newcommand{\calX}{{\cal X}} 
\newcommand{\calXX}{{\cal X\mbox{\raisebox{.3ex}{$\!\!\!\!\!-$}}}} 
\newcommand{\calXXX}{{\cal X\!\!\!\!\!-}} 
\newcommand{\gi}{{\raisebox{.0ex}{$\scriptscriptstyle{\cal X}$} 
\raisebox{.1ex} {$\scriptstyle{\!\!\!\!-}\:$}}} 
\newcommand{\intsim}{\int_0^1\!\!\!\!\!\!\!\!\!\sim} 
\newcommand{\intsimt}{\int_0^t\!\!\!\!\!\!\!\!\!\sim} 
\newcommand{\pp}{{\partial}} 
\newcommand{\al}{{\alpha}} 
\newcommand{\sB}{{\cal B}} 
\newcommand{\sL}{{\cal L}} 
\newcommand{\sF}{{\cal F}} 
\newcommand{\sE}{{\cal E}} 
\newcommand{\sX}{{\cal X}} 
\newcommand{\R}{{\rm I\!R}} 
\renewcommand{\L}{{\rm I\!L}} 
\newcommand{\vp}{\varphi} 
\newcommand{\N}{{\rm I\!N}} 
\def\ooo{\lip} 
\def\ccc{\rip} 
\newcommand{\ot}{\hat\otimes} 
\newcommand{\rP}{{\Bbb P}} 
\newcommand{\bfcdot}{{\mbox{\boldmath$\cdot$}}} 
 
\renewcommand{\varrho}{{\ell}} 
\newcommand{\dett}{{\textstyle{\det_2}}} 
\newcommand{\sign}{{\mbox{\rm sign}}} 
\newcommand{\TE}{{\rm TE}} 
\newcommand{\TA}{{\rm TA}} 
\newcommand{\E}{{\rm E\,}} 
\newcommand{\won}{{\mbox{\bf 1}}} 
\newcommand{\Lebn}{{\rm Leb}_n} 
\newcommand{\Prob}{{\rm Prob\,}} 
\newcommand{\sinc}{{\rm sinc\,}} 
\newcommand{\ctg}{{\rm ctg\,}} 
\newcommand{\loc}{{\rm loc}} 
\newcommand{\trace}{{\,\,\rm trace\,\,}} 
\newcommand{\Dom}{{\rm Dom}} 
\newcommand{\ifff}{\mbox{\ if and only if\ }} 
\newcommand{\nproof}{\noindent {\bf Proof:\ }} 
\newcommand{\remark}{\noindent {\bf Remark:\ }} 
\newcommand{\remarks}{\noindent {\bf Remarks:\ }} 
\newcommand{\note}{\noindent {\bf Note:\ }}

\newcommand{\boldx}{{\bf x}} 
\newcommand{\boldX}{{\bf X}} 
\newcommand{\boldy}{{\bf y}} 
\newcommand{\boldR}{{\bf R}} 
\newcommand{\uux}{\uu{x}} 
\newcommand{\uuY}{\uu{Y}} 
 
\newcommand{\limn}{\lim_{n \rightarrow \infty}} 
\newcommand{\limN}{\lim_{N \rightarrow \infty}} 
\newcommand{\limr}{\lim_{r \rightarrow \infty}} 
\newcommand{\limd}{\lim_{\delta \rightarrow \infty}} 
\newcommand{\limM}{\lim_{M \rightarrow \infty}} 
\newcommand{\limsupn}{\limsup_{n \rightarrow \infty}} 
 
\newcommand{\ra}{ \rightarrow }

\newcommand{\ARROW}[1] 
  {\begin{array}[t]{c}  \longrightarrow \\[-0.2cm] \textstyle{#1} \end{array} } 
 
\newcommand{\AR} 
 {\begin{array}[t]{c} 
  \longrightarrow \\[-0.3cm] 
  \scriptstyle {n\rightarrow \infty} 
  \end{array}} 
 
\newcommand{\pile}[2] 
  {\left( \begin{array}{c}  {#1}\\[-0.2cm] {#2} \end{array} \right) } 
 
\newcommand{\floor}[1]{\left\lfloor #1 \right\rfloor} 
 
\newcommand{\mmbox}[1]{\mbox{\scriptsize{#1}}} 
 
\newcommand{\ffrac}[2] 
  {\left( \frac{#1}{#2} \right)} 
 
\newcommand{\one}{\frac{1}{n}\:} 
\newcommand{\half}{\frac{1}{2}\:} 
 
\def\le{\leq} 
\def\ge{\geq} 
\def\lt{<} 
\def\gt{>} 
 
\def\squarebox#1{\hbox to #1{\hfill\vbox to #1{\vfill}}} 
\newcommand{\nqed}{\hspace*{\fill} 
           \vbox{\hrule\hbox{\vrule\squarebox{.667em}\vrule}\hrule}\bigskip} 
 
\title{Martingale representation for degenerate diffusions}

\author{ A. S. \"Ust\"unel} 
\maketitle 
\noindent 
{\bf Abstract:}{\small{Let $(W,H,\mu)$ be the classical Wiener space
    on $\R^d$.  Assume that $X=(X_t)$ is a diffusion process satisfying
    the stochastic differential equation $dX_t=\sigma(t,X)dB_t+b(t,X)dt$, where
  $\sigma:[0,1]\times C([0,1],\R^n)\to \R^n\otimes \R^d$,
  $b:[0,1]\times C([0,1],\R^n)\to \R^n$, $B$ is an $\R^d$-valued
  Brownian motion. We suppose that the weak
  uniqueness of this equation holds for any initial condition. We
  prove that any square integrable martingale $M$ w.r.t. to the filtration
  $(\calF_t(X),t\in [0,1])$ can be represented as 
$$
M_t=E[M_0]+\int_0^t P_s(X)\alpha_s(X).dB_s
$$
where $\alpha(X)$ is an $\R^d$-valued process adapted to
$(\calF_t(X),t\in [0,1])$, satisfying
$E\int_0^t(a(X_s)\alpha_s(X),\alpha_s(X))ds<\infty$, $a=\sigma^\star\sigma$
and $P_s(X)$ denotes a measurable version of the orthogonal projection
from $\R^d$ to $\sigma(X_s)^\star(\R^n)$. In particular, for any $h\in
H$, we have
\begin{equation}
\label{wick}
E[\rho(\delta
h)|\calF_1(X)]=\exp\left(\int_0^1(P_s(X)\dot{h}_s,dB_s)-\half\int_0^1|P_s(X)\dot{h}_s|^2ds\right)\,,
\end{equation}
where $\rho(\delta h)=\exp(\int_0^1(\dot{h}_s,dB_s)-\half
|H|_H^2)$. In the case the process $X$ is adapted to the Brownian
filtration, this result gives a new development as an
infinite series of the $L^2$-functionals of the degenerate
diffusions. We also give an adequate  notion of ``innovation
process'' associated to a degenerate diffusion which corresponds to
the strong solution when the Brownian motion is replaced by an adapted
perturbation of identity. This latter result gives the solution
of the causal Monge-Amp\`ere equation.}
\vspace{0.5cm} 


\section{\bf{Introduction} }
The representation of random variables with the stochastic integrals
with respect to some basic processes has a long history and also it
has very important applications, e.g., in signal theory, filtering,
optimal control, finance, stochastic differential equations, in
physics, etc. Let us recall the question: assume that we are given a
certain semimartingale $X$ indexed by $[0,1]$ for example. Let
$\calF(X)=(\calF_t(X),t\in [0,1])$ denote its filtration. The question is under
which conditions can we represent any martingale adapted to the
filtration of $X$ as a stochastic integral w.r.t. a fixed martingale
of the filtration $\calF(X)$? For the case of Wiener process, this
question has a very long history and it is almost impossible to give
an exact account of the contributing works. To our
knowledge, it has  been answered for the first time in the work of K. It\^o,
cf. \cite{Ito,W}. In \cite{CD}, C. Dellacherie has given a different
point of view to prove the representation theorem for the Wiener and
Poisson processes, based on the uniqueness of their laws. The case of nondegenerate
diffusion processes has been elucidated in \cite{S-V}, also studied in
\cite{E-H,JJ} (cf. also the references there) with also some remarks
about the degenerate case. The general case, using the notion of
multiplicity is given in \cite{D-V}. 

Although the nondegenerate case
is completely settled, the degenerate case has not been solved in a 
definitive way, in the sense that a unique minimal martingale with
respect to which the above mentioned property of representation holds, 
has not been discovered. In this work we are doing exactly this: we
prove that for a degenerate diffusion whose law is unique, there
exists a minimal martingale with respect to which every square
integrable, $\calF_1(X)$-measurable functional can be represented as a
stochastic integral of an $\calF(X)$-adapted process. 
To do this we first prove the density  of a class of
$\calF_1(X)$-measurable  stochastic integrals in $L^2(\calF_1(X))$ using the method
launched by C. Dellacherie, then we show that a
sequence from this class which is approximating any element of $L^2(\calF_1(X))$
defines an $L^2$-converging  sequence of processes with values in the range of the
adjoint of the diffusion coefficient, i.e., with values in
$\sigma^\star(t,X)(\R^n)$. Because of the degeneracy of $\sigma$, we
can not determine the limit of this sequence but its image under the
orthogonal projection $P_t(X): \R^d\to \sigma^\star(t,X)(\R^n)$ and
its limit is perfectly well-defined and in this way we see that the
minimal martingale for the representation of the functionals of the
diffusion (It\^o process) is nothing but $dm_t=P_t(X)dB_s$ in its
infinitesimal It\^o form, where $B$ is the Brownian motion governing
the process $X$. Note that even the adaptability of $(m_t,t\in [0,1])$
to the filtration  $\calF(X)$ is not evident. This result of
representation gives also existence of non-orthogonal chaos
representation of the elements of $L^2(\calF_1(X))$ as the multiple
ordered integrals with respect to the martingale $(m_t,t\in [0,1])$ of elements of $\{L^2(C_n,dt_1\times\cdots\times
dt_n),n\geq 1\}$, where
$C_n=\{(t_1,\ldots,t_n)\in[0,1]^n:t_1>\ldots>t_n\}$. 

In the third section we
define the notion of innovation process associated to a degenerate
diffusion process\footnote{Note that we use the word diffusion in the
  large sense, i.e., without demanding a Markov property.}}. In this
case the definition of innovation process is different from the
classical case of the perturbation of identity since we have to take
into account also the action of the projection operator-valued process
$(P_t(X), \,t\in [0,1])$. In Section 3 we extend the innovation
representation theorem of Fujisaki-Kallianpur-Kunita, \cite{FKK}, under the
hypothesis of strong exitence and uniqueness to
the case of degenerate diffusions. This result confirms the validity of
the choice that we have done to define the innovation process. In
particular, using this innovation process we can calculate the
conditional expectation of the Girsanov exponential of an adapted
drift $u$ with respect to the sigma algebra $\calF_1(X^U)$, where
$X^U$ is the solution of the diffusion stochastic differential
equation where the random input is equal to $U=B+u$. 

The results of the third section is then applied to the solution of
the adapted Monge-Amp\`ere equation in the case of the degenerate
diffusions, which extends the results of
\cite{FUZ,ASU-3,ASU-4}, cf. also \cite{fandu1,fandu2,fandu3}. In particular we
calculate the relative entropy of the law of $X^U$ with respect to the
law of $X$ by the use of preceding results. 

Let us note to finish this introduction that most of these results are
easily extendible to more general situations, for example the strong
existence and uniqueness hypothesis can be weakened in the entropic
calculations, we have tried to follow maximum homogeneity in the
hypothesis and these possible extensions may be treated in seperate
works.

\section{\bf{Stochastic integral representation of  functionals of
    diffusions }}
\noindent
Let $X=(X_t,t\in [0,1])$ be a weak  solution of the following
stochastic differential equation:
\begin{equation}
\label{SDE}
dX_t=b(t,X)dt+\sigma(t,X)dB_t,\,X_0=x,
\end{equation}
where $B=(B_t,t\in [0,1])$ is an $\R^d$-valued  Brownian motion and
$\sigma: [0,1]\times C([0,1],\R^n)\to L(\R^d,\R^n)$ and $b:
[0,1]\times C([0,1],\R^n)\to \R^n$ are  measurable maps, adapted to the
natural filtration of $C([0,1],\R^n)$ and of linear
growth. Recall that the classical theorem of Yamada-Watanabe says that the strong
uniqueness implies the uniqueness in law of the above SDE,
cf.\cite{Y-W,I-W}. Hence weak uniqueness is easier to obtain in the
applications. In this section we shall assume the weak
uniqueness of the equation and prove the martingale representation
property for the case where $\sigma$ may be degenerate.

More precisely, let $(\calF_t(X),t\in [0,1])$
be the filtration of $X$ and let us denote by $K$ the set of
$\R^n$-valued, $(\calF_t(X),t\in [0,1])$-adapted processes $\alpha(X)$,
s.t.
$$
E\int_0^1(a(s,X)\al_s(X),\al_s(X))ds<\infty\,,
$$
where $a(s,w)=\sigma(s,w)\sigma^\star(s,w),\,s\in[0,1],\,w\in C([0,1],\R^n)$. 
\begin{theorem}
\label{rep-thm}
The set $\Gamma=\{N\in
L^2(\calF_1(X)):\,N=E[N]+\int_0^1(\al_s(X),\sigma(s,X)dB_s),\,\al\in K\}$ is
dense in $L^2(\calF_1(X))$.
\end{theorem}
\nproof Suppose that there is some $M\in L^2(\calF_1(X))$  which is orthogonal (in
$L^2$) to $\Gamma$. Using the
usual stopping technique, we can assume that the corresponding
$(\calF_t(X),t\in [0,1])$-martingale is bounded and  positive whose expectation is equal to
one. The orthogonality implies that
$(M_t(f(X_t)-\int_0^tLf(s,X)ds),t\in[0,1])$ is again a (local) martingale for
any smooth function $f:\R^n\to\R$, where $L$ is defined as 
\begin{equation}
\label{generator}
Lf(t,X)=\half\sum_{i,j}a_{i,j}(t,X)\partial_{i,j}f(X_t)+\sum_i
b_i(t,X)\partial_i f(X_t)\,.
\end{equation}
This implies, with Proposition IV.2.1 of \cite{I-W}, that, under the
measure $M\cdot P$, $X$ is again a weak solution of \ref{SDE}. 
By the uniqueness in law, we get $X(M\cdot P)=X(P)$, since $M$ is
$\calF_1(X)$-measurable, we should have $M=1$. Hence the functionals
of the diffusion which are orthogonal to the above set of stochastic
integrals are almost surely constant. Consequently, the set $\Gamma$
is total in $L^2(\calF_1(X))$. 
\nqed

\noindent
The following is the {\bf{extension of the martingale representation
theorem to the functionals of degenerate diffusions}}:
\begin{theorem}
\label{adapt_thm}
Denote by $P_s(X)$ a measurable version of the orthogonal projection
from $\R^d$ onto $\sigma(s,X)^\star(\R^n)\subset\R^d$ and let $F\in
L^2(\calF_1(X))$ be any random variable with zero expectation. Then
there exists a process $\xi(X)\in L_a^2(dt\times dP;\R^d)$, adapted to
$(\calF_t(X),t\in [0,1])$, such that
$$
F(X)=\int_0^1 (P_s(X)\xi_s(X),dB_s)_{\R^d}=\int_0^1 (\xi_s(X),P_s(X)
dB_s)_{\R^n}
$$
a.s.\\
Conversely, any stochastic integral of the form
$$
\int_0^1(P_s(X)\xi_s(X),dB_s)\,,
$$
where $\xi(X)$ is an $(\calF_t(X),t\in [0,1])$-adapted, measurable
process with $E\int_0^1|P_s(X)\xi_s(X)|^2ds<\infty$, gives rise to an
$\calF_1(X)$-measurable random variable.
\end{theorem}
\nproof
From Theorem \ref{rep-thm}, there exists a sequence $(F_n(X),n\geq
1)\subset \Gamma$, where
$\Gamma$ is defined in the statement of Theorem \ref{rep-thm},
 converging to $F(X)$ in $L^2$. We can suppose that $E[F_n(X)]=0$, for
 any $n\geq 1$. Hence
$$
F_n(X)= \int_0^1(\ga^n_s(X),\sigma(s,X) dB_s)_{\R^d}=\int_0^1(\sigma^\star(s,X)\ga^n_s(X),dB_s)_{\R^d}\,.
$$
As explained above $\ga^n$ is an $(\calF_t(X),t\in [0,1])$-adapted, $\R^n$-valued
process satisfying
$E\int_0^1(a(X_s)\ga^n_s(X),\ga^n_s(X))ds<\infty$. 
Since $F_n(X)\to F(X)$ in $L^2$, 
$$
\lim_{n,m\to\infty}E\int_0^1|\sigma^\star(s,X)\ga^n_s(X)-\sigma^\star(s,X)\ga^m_s(X)|^2ds=0\,.
$$
Let $\alpha^n_s(X)=\sigma^\star(X_s)\ga^n_s(X)$, as
$P_s(X)\alpha^n_s(S)=P_s(X)
\sigma^\star(X_s)\ga^n_s(X)=\sigma^\star(X_s)\ga^n_s(X)$, 
$(P_s(X)\alpha^n_s(X),n\geq 1)$ converges to some $\xi_s(X)$ in
$L_a^2(ds\times dP;\R^d)$. As $P_s(X)$ is an orthogonal projection,
$(P_s(X)\alpha_s^n(X),n\geq 1)$ converges to $P_s(X)\xi_s(X)$ also in
$L_a^2(ds\times dP;\R^d)$. Therefore
\beaa
\int_0^1P_s(X)\xi_s(X).dB_s&=&\lim_n\int_0^1P_s(X)\alpha^n_s(X).dB_s\\
&=&\lim_n\int_0^1(\ga^n_s(X),\sigma(s,X)dB_s)\\
&=&\lim_nF_n(X)=F(X)\,.
\eeaa
Let now $G\in L^2(P)$ be given by $G=\int_0^1(P_s(X)\eta_s(X),dB_s)$
and assume that it is not $\calF_1(X)$-measurable. Then $G-E[G|\calF_1(X)]$ is
orthogonal to $L^2(\calF_1(X))$. It follows from the first part
of the theorem that we can represent $E[G|\calF_1(X)]$ as
$\int_0^1(P_s(X)\xi_s(X),dB_s)$. Let us define $h=\eta-\xi$, then the
orthogonality mentioned above implies that 
$$
E\left[\int_0^1(P_s(X)h_s,dB_s).\int_0^1(\alpha_s(X),\sigma(s,X)dB_s)\right]=0
$$
for any $(\calF_t(X),t\in [0,1])$-adapted, measurable $\alpha$ such
that $E\int_0^1(a(s,X)\alpha_s,\alpha_s)ds<\infty$. Consequently
$P_s(X)h_s=0$ $ds\times dP$-a.s., hence $G=E[G|\calF_1(X)]$ P-a.s.

\nqed

\remark Let $\eta$ be an adapted process such that $\eta_s$ belongs to
the orthogonal complement of $\sigma^\star(\R^n)$ in $\R^d$ $ds\times
dP$-a.s. Then $\eta+\xi$  can also be used to represent
$F(X)$. Hence $\xi(X)$ is not unique but $P(X)\xi(X)$ is always unique.

\begin{theorem}
\label{iden-thm}
Let $\dot{u}\in L^2(dt\times dP,\R^d)$ be adapted to the Brownian filtration, then we have
\begin{equation}
\label{cond-exp}
E\left[\int_0^1(\dot{u}_s,dB_s)|\calF_1(X)\right]=\int_0^1(E[P_s(X)\dot{u}_s|\calF_s(X)],dB_s)
\end{equation}
almost surely.
\end{theorem}
\nproof
We have to prove first  that the right hand side of (\ref{cond-exp})
is $\calF_1(X)$-measurable. We know from Theorem \ref{adapt_thm} that
the left side of (\ref{cond-exp}) can be represented as 
$$
E\left[\int_0^1(\dot{u}_s,dB_s)|\calF_1(X)\right]=\int_0^1 P_s(X)\xi_s(X).dB_s
$$
for some $\xi\in L_a^2(dt\times dP;\R^d)$. Let
$F(X)=\int_0^1P_s(X)\alpha_s(X).dB_s$ be any element of
$L^2(\calF_1(X))$ with $\alpha(X)\in L_a^2(dt\times dP;\R^d)$. We have
\beaa
E\left[\left(\int_0^1\dot{u}_s.dB_s\right) F(X)\right]&=&E\int_0^1(\dot{u}_s,P_s(X)\alpha_s(X))ds\\
&=&E\int_0^1(P_s(X)E[\dot{u}_s|\calF_s(X)],P_s(X)\alpha_s(X))ds\\
&=&E\left[\left(\int_0^1 P_s(X)\xi_s(X).dB_s\right)\,F(X)\right]\\
&=&E\int_0^1 (P_s(X)\xi_s(X),P_s(X)\alpha_s(X))ds\,,
\eeaa
hence $P_s(X)\xi_s(X)=P_s(X)E[\dot{u}_s|\calF_s(X)]$ $ds\times
dP$-a.s., in particular 
$$
\int_0^1(E[P_s(X)\dot{u}_s|\calF_s(X)],dB_s)
$$
is $\calF_1(X)$-measurable. The above identification assures then the
validity of the relation (\ref{cond-exp}).
\nqed

\begin{corollary}
\label{cor_1}
Let $h\in H^1([0,1],\R^d)$ (i.e., the Cameron-Martin space), denote by
$\rho(\delta h)$ the Wick exponential
$\exp(\int_0^1(\dot{h}_s,dB_s)-\half\int_0^1|\dot{h}_s|^2ds)$, 
 then we have
$$
E[\rho(\delta h)|\calF_1(X)]=
\exp\left(\int_0^1(P_s(X)\dot{h}_s,dB_s)-\half\int_0^1|P_s(X)\dot{h}_s|^2ds\right)\,.
$$
\end{corollary}

\begin{theorem}
\label{chaos-thm}
\begin{itemize}
\item Assume that $\calF_t(X)\subset \calF_t(B)$ for any $t\in [0,1]$, where
$(\calF_t(B),t\in [0,1])$ represents the filtration of the Brownian motion. Define the martingale $m=(m_t,t\in [0,1])$ as
$m_t=\int_0^tP_s(X) dB_s$,
then the set 
$$
K=\{\rho(\delta_m(h)):\,h\in H\}
$$ 
is total in
$L^2(\calF_1(X))$, where
$\rho(\delta_m(h))=\exp\left(\int_0^1(\dot{h}_s,dm_s)-\half\int_0^1|P_s(X)\dot{h}_s|^2ds\right)\,$.
In particular, any element $F$ of $L^2(\calF_1(X))$ can be written in
a unique way as the sum
\begin{equation}
\label{m_wi}
F=E[F]+\sum_{n=1}^\infty\int_{C_n}(f_n(s_1,\ldots,s_n),dm_{s_1}\otimes\ldots\otimes
dm_{s_n})
\end{equation}
where $C_n$ is the $n$-dimensional simplex in $[0,1]^n$ and $f_n\in
L^2(C_n, ds^{\otimes n})\otimes (\R^d)^{\otimes n}$.
\item   
More generally, without the hypothesis $\calF_t(X)\subset \calF_t(B)$,
for any $F\in L^2(\calF_1(X))\cap L^2(\calF_1(B))$, the conclusions of
the first part of the theorem hold true.
\end{itemize}
\end{theorem}
\nproof
Let $F\in L^2(\calF_1(X))$, assume that $F$ is orthogonal to $K$,
i.e. $E[F\,\rho(\delta_m(h))]=0$
for any $h\in H$. From Corollary \ref{cor_1}  $\rho(\delta_m(h))=E[\rho(\delta h)|\calF_1(X)]$,
hence
$$
E[F\,\rho(\delta h)]=E[F\,\rho(\delta_m(h))]=0\,,
$$
hence $F$ is also orthogonal to $\mathcal{E}=\{\rho(\delta h):\,h\in H\}$,
which is total in $L^2(\calF_1(B)$, where $\calF_1(B)$ is the
$\sigma$-algebra generated by the governing Brownian motion, therefore
$F=0$. Consequently the span of $K$ is dense in $L^2(\calF_1(X))$.
Theorem \ref{iden-thm} and Corollary \ref{cor_1} allow us to calculate the conditional
expectations of the multiple Wiener integrals w.r.t. $\calF_1(X)$ and
the result will be multiple iterated stochastic integrals w.r.t. $m$
in the form given by the formula (\ref{m_wi});
here we have to be careful as the operator valued process
$(P_s(X),s\in [0,1])$ is not deterministic, the symmetric
interpretation of the Ito-Wiener integrals w.r.t. the Brownian motion
is no longer valid in our case and they have to be written as iterated
Ito integrals.
 \nqed

\begin{remarkk}
In this theorem we need to make the hypothesis $\calF_t(X)\subset
\calF_t(B)$ for any $t\in [0,1]$ to assure the non-symmetric chaos
representation (\ref{m_wi}). Without this hypothesis, although we have
Theorem \ref{adapt_thm} and when we iterate it we have a similar
representation, but, consisting of a finite number of terms. It is
not possible to push this procedure up to infinity since we have no
control at infinity.
\end{remarkk}
\begin{remarkk}
Let us note that if there is no strong solution to the equation
defining the process $X$, the chaotic  representation property may
fail. For example, let $U$ be a weak solution of 
\begin{equation}
\label{co-ex}
dU_t=\alpha_t(U)dt+dB_t\,,
\end{equation}
with $U_0$ given. Assume that (\ref{co-ex}) has no strong solution, as
it may happen in the famous example of Tsirelson (\cite{I-W}), i.e.,
$U$ is not measurable w.r.t. the sigma algebra generated by $B$, then
we have no chaotic representation property for the elements of
$L^2(\calF_1(U))$ in terms of the iterated stochastic integrals of
deterministic functions on $C_n, n\in \N$ w.r.t. $B$; the contrary
would imply the equality of $\calF_1(U)$ and of $\calF_1(B)$, which
would contradict the non-existence of strong solutions. 
\end{remarkk}

\section{\bf{Innovation Process}}
\noindent
In this section {\bf{we assume that the SDE \ref{SDE} has unique strong
solution}}. To fix the ideas, we can suppose that $\sigma$ and $b$ satisfy the
following kind of Lipschitz condition on the path space: for
$\xi,\,\eta\in C([0,1],\R^n)$:
$$
|\ga(t,\xi)-\ga(t,\eta)|\leq K\sup_{s\leq t}|\xi(s)-\eta(s)|
$$
for $\ga$ being equal  either to $\sigma$ or to $b$ with corresponding
Euclidean norms at the left hand side.

\noindent
Assume that $\dot{u}\in L^2(dt\times dP,\R^d)$ is a process adapted to
the filtration of Brownian motion, let $U=(U_t,t\in [0,1])$ be defined
as $U_t=B_t+\int_0^t\dot{u}_sds$. We denote by $X^U$ the strong
solution of the equation 
\beaa
dX^U_t&=&\sigma(t,X^U) dU_t+b(t,X^U)dt\\
&=&\sigma(t,X^U)(dB_t+\dot{u}_tdt)+b(t,X^U)dt
\eeaa 
Since $B$ is the canonical Brownian motion, we
have $X_t^U=X_t\circ U$ a.s. Moreover, if $\eta\in \R^n$ is any
vector, we have
\beaa
(\sigma(t,X^U) \dot{u}_t,\eta)&=&(\dot{u}_t, \sigma(t,X^U)
^\star\eta)\\
&=&(P_t(X^U)\dot{u}_t, \sigma(t,X^U)^\star\eta)\\
&=&(\sigma(t,X^U)P_t(X^U)\dot{u}_t, \eta)\,,
\eeaa
where $P_t(X^U)$ denotes the orthogonal projection from $\R^d$ onto
$\sigma(t,X^U)^\star(\R^n)$. Hence, in order to define a reasonably
useful concept of innovation, we should estimate the perturbation
$\dot{u}$ simultaneously w.r.t. the both projections, i.e., with
respect to the conditional expectation $E[\cdot|\calF_s(X^U)]$ (which
is a projection) and also w.r.t. $P_s(X^U)$: Let $Z=(Z_t,t\in[0,1])$
(to avoid the ambiguity we shall also use the notation $Z^U$ if necessary)
be defined as
\begin{equation}
\label{innov}
Z_t=B_t+\int_0^t(\dot{u}_s-E[P_s(X^U)\dot{u}_s|\calF_s(X^U)])ds\,,
\end{equation}
where $(\calF_s(X^U),s\in [0,1])$ denotes the filtration of $X^U$. 
\begin{proposition}
\label{inno_adapt}
The process $(\int_0^tP_s(X^U)dZ_s,t\in [0,1])$ is an $(\calF_t(X^U),t\in [0,1])$-local martingale.
\end{proposition}
\nproof
Assume to begin that $|u|_H^2=\int_0^1|\dot{u}_s|^2ds\in
L^\infty(P)$. Let us first prove that the process under question is adapted to the
filtration $(\calF_t(X^U),t\in [0,1])$:
From the Girsanov theorem, the process $(\int_0^t\sigma (s,X^U)dU_s,t\in
[0,1])$ is adapted to the filtration  $(\calF_t(X^U),t\in [0,1])$,
using the same method as in Theorem \ref{adapt_thm} by replacing $B$
by $U$ and the probability $dP$ by $\rho(-\delta u) dP$, we conclude
that the process $(\int_0^tP_s(X^U)dU_s,t\in [0,1])$, and hence the
process $(\int_0^tP_s(X^U)dZ_s,t\in [0,1])$ is adapted to the
filtration  $(\calF_t(X^U),t\in [0,1])$. To show the (local)
martingale property it suffices to write that 
$$
\int_0^tP_s(X^U)dZ_s=\int_0^tP_s(X^U)dB_s+\int_0^tP_s(X^U)\left[\dot{u}_s-E[\dot{u}_s|\calF_s(X^U)]\right]ds
$$
from which the martingale property follows. The general case follows
from a stopping argument.

\nqed

\remark
Note that $Z=Z^U$ is not a Brownian motion.
\begin{theorem}
\label{cond_exp}
Let $\dot{u}\in L^2(dt\times dP,\R^d)$ be such that $E[\rho(-\delta
u)]=1$, then we have
\beaa
\zeta_t&=&E[\rho(-\delta u)|\calF_t(X^U)]\\
&=&\exp\left(-\int_0^tP_s(X^U) E[\dot{u}_s|\calF_s(X^U)]\cdot
  dZ_s-\half\int_0^t|P_s(X^U)E[\dot{u}_s|\calF_s(X^U)]|^2ds\right)\,.
\eeaa
\end{theorem}
\nproof
Assume first that $|u|_H^2=\int_0^1|\dot{u}_s|^2ds\in
L^\infty(P)$. Let $(X^U_t,t\in [0,1])$ be the (strong) solution of $dX_t=\sigma
(X_t)dU_t$, where $dU_t=dB_t+\dot{u}_tdt$ and 
let $f$ be a $C^2$-function on $\R^n$. Using
the It\^o formula, we calculate the Doob-Meyer process associated to
the semimartingales $(\zeta_t)$ and $(f(X^U_t))$:
\begin{equation}
\label{D_M}
\langle f\circ X^U,\zeta\rangle_t=-\int_0^t\zeta_s\left(Df(X^U_s),\sigma
(s,X^U)P_s(X^U) E[\dot{u}_s|\calF_s(X^U)] \right)ds\,.
\end{equation}
Let $f\in C^2$, using again the It\^o formula and the relation (\ref{D_M}) we get
\beaa
f(X^U_t)\zeta_t&=&\int_0^tf(X_s^U)d\zeta_s+\int_0^t\zeta_s
(Lf)(s,X^U)ds+\int_0^t\zeta_s(Df(X^U_s),\sigma(s,X^U)dU_s)\\
&&-\int_0^t\zeta_s (Df(X^U_s),\sigma(s,X^U)P_s(X^U)E[\dot{u}_s|\calF_s(X^U)])ds\,,
\eeaa
where $Lf$ is defined by the relation \ref{generator}. Therefore  the process 
$$
\left(f(X^U_t)\zeta_t-\int_0^t (Lf)(s,X^U))\zeta_sds,t\in [0,1]\right)
$$
is a $P$-local martingale, therefore,  by the uniqueness in law of the solution, we
should have
$$
E[\zeta_1 F(X^U)=E[F(X)]=E[F(X^U)\rho(-\delta u)]
$$ 
for any  $F\in C_b$ on the path space, and  the last inequality follows from the Girsanov theorem. It
suffices then to remark that $\zeta_1$ is
$\calF_1(X^U)$-measurable by Theorem \ref{iden-thm} and again by the
Girsanov theorem which allows us to replace $B$ by $U$. The general case follows from the usual
stopping time argument: let $T_n=\inf(t: \int_0^t|\dot{u}_s|^2ds\geq
n)$ and define $\dot{u}_s^n=1_{[0,T_n]}(s)\dot{u}_s$ and let
$U^n=B+\int_0^\cdot \dot{u}_s^nds$. Then $X^{U^n}$ converges almost
surely uniformly to $X^U$, hence
$\lim_nE[\,\cdot\,|\calF_t(X^{U^n})]=E[\,\cdot\,|\calF_t(X^{U})]$ $dt$-almost
surely as bounded operators on $L^1(P)$ and $(P_s(X^{U^n}),n\geq 1)$
converges to $P_s(XÛ)$ $ds\times dP$-almost surely.  Moreover $(\rho(-\delta u^n),n\geq
1)$ converges strongly in $L^1(P)$ to $\rho(-\delta u)$, therefore the general
case follows.
\nqed

\noindent
The following result is the generalization of the celebrated
innovation's theorem to the degenerate case, cf.\cite{FKK}:
\begin{theorem}
\label{innovation-thm}
Let $(M_t,t\in [0,1])$ be a square integrable $(P,(\calF_t(X^U),t\in
[0,1]))$-martingale, then it can be represented as a stochastic
integral of an $(\calF_t(X^U),t\in [0,1])$-adapted, $\R^d$-valued process
$\beta(X^U)$ in the following way:
$$
M_t=M_0+\int_0^t (P_s(X^U)\beta_s(X^U),dZ_s)
$$
$P$-a.s, where 
$$
\int_0^t|P_s(X^U)\beta_s(X^U)|^2ds<\infty
$$
$P$-a.s., for any $t\in [0,1]$.
\end{theorem}
\nproof
Assume that $M$ is a  $(P,(\calF_t(X^U),t\in[0,1]))$-martingale, then
for any $s<t$ and $A\in \calF_s(X^U)$ we have
\beaa
E\left[\frac{M_t}{\zeta_t}1_A\rho(-\delta u)\right]&=&E
\left[\frac{M_t}{\zeta_t}1_A\zeta_t\right]\\
&=&E[M_t1_A]=E[M_s1_A]\\
&=&E\left[\frac{M_s}{\zeta_s}1_A\rho(-\delta u)\right]
\eeaa
where $\zeta$ is the optional projection of $\rho(-\delta u)$
w.r.t. the filtration $(\calF_t(X^U),t\in[0,1])$ as calculated in
Theorem \ref{cond_exp}. Consequently $(M_t/\zeta_t,t\in [0,1])$ is a
$(Q,(\calF_t(X^U),t\in[0,1]))$-martingale, where $dQ=\rho(-\delta
u)dP$. As $U$ is a $Q$-Brownian motion, from Theorem \ref{adapt_thm},
we can represent $(M_t/\zeta_t,t\in [0,1])$ as
$$
\frac{M_t}{\zeta_t}=c+\int_0^t(P_s(X^U)\dot{\alpha}_s(X^U),dU_s)\,,
$$
then using the It\^o formula 
\beaa
M_t&=&\frac{M_t}{\zeta_t}\zeta_t\\
&=&c+\int_0^t\zeta_s(P_s(X^U)\dot{\alpha}_s(X^U),dU_s)-\int_0^t\frac{M_s}{\zeta_s}\zeta_s(P_s(X^U) E[\dot{u}_s|\calF_s(X^U)],dZ_s)\\
&&-\int_0^t\zeta_s(P_s(X^U)\dot{\alpha}_s(X^U),E[\dot{u}_s|\calF_s(X^U)])ds\\
&=&c+\int_0^t\zeta_s(P_s(X^U)\dot{\alpha}_s(X^U),dZ_s+P_s(X^U)
E[\dot{u}_s|\calF_s(X^U)])ds-\int_0^tM_s (P_s(X^U)
E[\dot{u}_s|\calF_s(X^U)],dZ_s)\\
&&-\int_0^t\zeta_s(P_s(X^U)\dot{\alpha}_s(X^U),
E[\dot{u}_s|\calF_s(X^U)])ds\\
&=&c+\int_0^t\left(P_s(X^U)\left[\zeta_s\dot{\alpha}_s(X^U)-M_s
E[\dot{u}_s|\calF_s(X^U)]\right],dZ_s\right)
\eeaa
and this completes the proof.
\nqed

\section{\bf{ Entropy Calculation and Monge-Amp\`ere Equation}}
\noindent
Assume that $l(X)$ is a probability density measurable
w.r.t. $\calF_1(X)$, i.e., $E[l(X)]=1$ and with finite entropy: 
$E[l(X)\log l(X)]<\infty$. We want to find a process $U=B+u=B+\int_0^\cdot \dot{u}_sds$ which  is an adapted
perturbation of the Brownian motion $B$ such that 
$$
l(X)=\frac{dX^U(P)}{dX(P)}\circ X\,.
$$
This problem is called the causal Monge-Amp\`ere problem. To simplify the calculations, we shall assume that $l\circ X$ is
$P$-a.s. strictly positive.  Assume that
such a $U$ (hence $u$) exists and that $u$ satisfies the Girsanov
theorem, i.e., $E[\rho(-\delta_Bu)]=1$. Then the Girsanov theorem
implies that 
\begin{equation}
\label{M-A.eqn}
l\circ X^U\,E[\rho(-\delta u)|\calF_1(X^U)]=1
\end{equation}
$P$-a.s., which is the causal version of the Monge-Amp\`ere
equation. From Theorem \ref{adapt_thm}, $l\circ X$ can be represented
as 
$$
l\circ X=\exp\left(-\int_0^1P_s(X)\dot{v}_s(X)\cdot dB_s-\half
  \int_0^1|P_s(X)\dot{v}_s(X)|^2ds\right)\,,
$$
where $(P_s(X)\dot{v}_s(X),s\in [0,1])$ is adapted to the filtration of $X$ and
$\int_0^1|P_s(X)\dot{v}_s(X)|^2ds<\infty$ $P$-a.s.
Besides, since $U$ is a Brownian motion under the
probability $\rho(-\delta u)\,dP$, it follows  from Theorem
\ref{adapt_thm} that  $l\circ X^U$ can be represented as 
\begin{equation}
\label{rep_formula}
l\circ X^U=\exp\left(-\int_0^1P_s(X^U)\dot{v}_s(X^U)\cdot dU_s-\half
  \int_0^1|P_s(X^U)\dot{v}_s(X^U)|^2ds\right)\,.
\end{equation}


\noindent
Inserting the right hand side of (\ref{rep_formula}) and $E[\rho(-\delta u)|\calF_1(X^U)]$ which
is already calculated in Theorem \ref{cond_exp} in the Monge-Amp\`ere
equation ({\ref{M-A.eqn}) and then taking the
logarithm of the final expression, we obtain 
\beaa
&&\int_0^1P_s(X^U)\left(E[\dot{u}_s|\calF_s(X^U)]+\dot{v}_s(X^U)\right)\cdot dZ_s\\
&&+\half\int_0^1\left(|P_s(X^U)\dot{v}_s(X^U)+P_s(X^U)
E[\dot{u}_s|\calF_s(X^U)]|^2\right)ds=0\,.
\eeaa
This relation implies that
\begin{equation}
\label{rep.eqn}
P_s(X^U_s)\left(\dot{v}_s(X^U)+E[\dot{u}_s|\calF_s(X^U)]\right)=0
\end{equation}
$ds\times dP$-almost surely,  which is a quite  elaborate nonlinear equation. From the Monge-Amp\`ere equation (\ref{M-A.eqn}) we can calculate the
relative entropy between $X^U(P)$ and $X(P)$, denoted by
$H(X^U(P)|X(P))$:

\begin{theorem}
\label{entropy}
Suppose that $l$ is an  $X(P)$-almost surely strictly positive
density. There exists some $\dot{u}\in L^2(ds\times dP)$ with $E[\rho(-\delta
u)]=1$ with 
$$
\frac{dX^U(P)}{dX(P)}=l
$$
if and only if
$$
P_s(X^U_s)\left(\dot{v}_s(X^U)+E[\dot{u}_s|\calF_s(X^U)]\right)=0.
$$
In this case we also have 
$$
H(X^U(P)|X(P))=\half
E\int_0^1|P_s(X^U)E[\dot{u}_s|\calF_s(X^U)]|^2ds=\half
E\int_0^1|P_s(X^U)\dot{v}_s(X^U)|^2ds\,.
$$
\end{theorem}
\nproof
\beaa
H(X^U(P)|X(P))&=&\int\log\frac{dX^U(P)}{dX(P)} dX^U(P)\\
&=&\int\log\frac{dX^U(P)}{dX(P)}\circ X^U dP\\
&=&\int\log l\circ X^U dP\\
&=&\half E\int_0^1|P_s(X^U)E[\dot{u}_s|\calF_s(X^U)]|^2ds\,,
\eeaa
provided that $\dot{u}\in L^2(ds\times dP)$ and the first equality
follows, the second one is a consequence of the relation
(\ref{rep.eqn}), it can be also proven directly from the Girsanov
theorem.
\nqed

\begin{proposition}
\label{SDE_monge}
Assume that $l$ and $u$ are given as above. Suppose furthermore that 
\begin{equation}
\label{hyp}
H(X^U(P)|X(P))=\half E\int_0^1|P_s(X^U)\dot{u}_s|^2ds\,.
\end{equation}
Then the following equation holds true:
\begin{equation}
\label{eqn_1}
P_s(X^U)dU_s+P_s(X^U)\dot{v}_s\circ X^Uds= P_s(X^U)dB_s
\end{equation}
almost surely. In particular, $X^U$ satisfies the following stochastic
differential equation:
\begin{equation}
\label{eqn_2}
dX^U_t=\sigma(t,X^U)(dB_t-\dot{v}_t\circ X^Udt)+b(t,X^U)dt,
\end{equation}
with the same initial condition as $X$.
\end{proposition}
\nproof
The hypothesis (\ref{hyp}) implies that the process
$(P_t(X^U)\dot{u}_t,t\in [0,1])$ is $ds$-almost surely adapted to the
filtration $(\calF_t(X^U),t\in [0,1])$, hence we get from the equality
(\ref{rep.eqn}) the relation
$$
P_t(X^U)(\dot{v}_t\circ X^U+\dot{u}_t)=0\,,
$$
which implies at once the relation (\ref{eqn_1}). To see the next one,
note that
\beaa
dX^U_t&=&\sigma(t,X^U)(dB_t+\dot{u}_tdt) +b(t,X^U)dt\\
&=&\sigma(t,X^U)(dB_t+P_t(X^U)\dot{u}_tdt) +b(t,X^U)dt\\
&=&\sigma(t,X^U)(dB_t-P_t(X^U)\dot{v}_t\circ X^Udt) +b(t,X^U)dt\\
&=&\sigma(t,X^U)(dB_t-\dot{v}_t\circ X^Udt) +b(t,X^U)dt
\eeaa
where we have used the fact that
$\sigma(t,X^U)\eta=\sigma(t,X^U)P_t(X^U)\eta$ for any vector in $\R^d$
since $P_t(X^U)$ is the orthogonal projection of $\R^d$ onto
$\sigma(X_t^U)^\star(\R^n)$.
\nqed

\noindent
Theorem \ref{entropy} can be extended as follows
\begin{theorem}
\label{ineq_entro}
Assume that $u\in L^2_a(dt\times dP,H)$ and denote by $U$ the process
$(B_t+\int_0^t \dot{u}_sds,t\in [0,1])$.  assume also, as before, the Lipschitz
hypothesis about the drift and diffusion coefficients, then the  following
inequality holds true:
\begin{equation}
\label{entro-ineq}
H(X^U(P)|X(P))\leq \half
E\int_0^1|P_s(X^U)E[\dot{u}_s|\calF_s(X^U)]|^2ds\,.
\end{equation}
\end{theorem}
\nproof
If $u\in L^\infty_a(dt\times dP,H)$, then the claim with equality
(instead of inequality) follows from Theorem \ref{entropy}. For the
case $u\in L^2_a(dt\times dP,H)$, define
$T_n=\inf(t>0:\int_0^t|\dot{u}_s|^2ds>n)$, then $u^n$ defined by
$$
u^n(t)=\int_0^t 1_{[0,T_n]}(s)\dot{u}_s ds
$$
is in $ L^\infty_a(dt\times dP,H)$, hence we have
\begin{equation}
\label{entro-eq}
H(X^{U^n}(P)|X(P))= \half
E\int_0^1|P_s(X^{U{^n}})E[1_{[0,T_n]}(s)\dot{u}_s|\calF_s(X^{U^n})]|^2ds\,.
\end{equation}
As $n\to \infty$, $(X^{U^n}(P),\,n\geq 1)$ converges weakly to $X^U(P)$
and the weak lower semi-continuity of the entropy implies that 
\beaa
H(X^U(P)|X(P))&\leq&\liminf_n\half
E\int_0^1|P_s(X^{U{^n}})E[1_{[0,T_n]}(s)\dot{u}_s|\calF_s(X^{U^n})]|^2ds\\
&=&\half E\int_0^1|P_s(X^U)E[\dot{u}_s|\calF_s(X^U)]|^2ds\,,
\eeaa
where the limit of the right hand side of the equation
(\ref{entro-eq}) follows from the Lipschitz hypothesis.
\nqed

\vspace{2cm}
\footnotesize{
\noindent
A. S. \"Ust\"unel, Bilkent University, Math. Dept., Ankara, Turkey\\
ustunel@fen.bilkent.edu.tr}


\begin{thebibliography}{99} 

\bibitem{D-V}
M. H. A. Davis and P. Varaiya: ``The multiplicity of an increasing
family of $\sigma$-fields''. The Annals of Probab.
{\bf{2}}, 958–963, 1974.


\bibitem{CD}
C. Dellacherie:`` Int´egrales stochastiques par rapport aux processus de Wiener 
ou de Poisson''. In S´eminaire de Probabilit´es VIII, Lecture Notes in Mathematics 
{\bf{381}}. Springer-Verlag, 1973. [Correction dans SP IX, LNM 465.] 





\bibitem{E-H}
N. El Karoui and H. Reinhard:
``Processus de diﬀusion dans $\R^n$''
S´eminaire de Probabilit´es, VII (Univ. Strasbourg, ann\'ee universitaire 
1971–-1972), pp. 95–117. Lecture Notes in Math., Vol. {\bf{321}}, Springer, 
Berlin, 1973. 




\bibitem{fandu1}
D. Feyel and A. S. \"Ust\"unel: ``Transport of measures on Wiener
space and the Girsanov theorem''. Comptes Rendus Math\'ematiques,
Vol. {\bf 334}, Issue 1, 1025-1028,  2002.


\bibitem{fandu2}
D. Feyel,  A.S. \"Ust\"unel:
Monge-Kantorovitch measure transportation and Monge-Amp\`ere equation on
Wiener space.
\emph{Probab. Theor. Relat. Fields}, {\bf 128}, no.~3, pp.~347--385,
2004.


\bibitem{fandu3} 
D. Feyel and A. S. \"Ust\"unel: ``Monge-Kantorovitch measure
transportation, Monge-Amp\`ere equation and the It\^o calculus''.  
 Advanced Studies in Pure Mathematics, Math. Soc. of Japan,
 Vol. {\bf{41}}, p. 32-49, 2004. Mathematical Society of Japan. 

\bibitem{FUZ}
D. Feyel, A.S. \"Ust\"unel and M. Zakai:
``Realization of Positive Random Variables via Absolutely Continuous
Transformations of Measure on Wiener Space''. Probability Surveys,Vol.
3, (electronic) p.170-205, 2006.

\bibitem{FKK}
M. Fujisaki, G. Kallianpur and H. Kunita:
``Stochastic Differential Equations for the Nonlinear Filtring
Problem''. Osaka J. Math, {\bf{9}}, p. 19-40, 1972.

\bibitem{I-W}
N. Ikeda and S. Watanabe: {\sl{Stochastic Differential Equations and
    Diffusion Processes}}


\bibitem{Ito}
K. It\^o:`` Multiple Wiener integrals''.  J. Math. Soc. Japan,
{\bf{3}}, p. 157–169, 1951.





\bibitem{JJ}
J. Jacod: {\sl{ Calcul stochastique et probl\'emes de martingales}}.  Lecture Notes in
Mathematics, Vol.{\bf714}. Springer-Verlag, 1979.



\bibitem{S-V}
D.W. Stroock and S.R.S. Varadhan:
``Diffusion processes with continuous coefficients 1''. Comm. Pure and
Appl. Math., {\bf{22}}, pp.345-400, 1969.





 \bibitem{ASU-3}
 A. S. \"Ust\"unel:``Entropy, invertibility and variational calculus
 of adapted shifts on Wiener space''. 
  J. Funct. Anal. {\bf{257}}, no. 8, p.3655--3689, 2009.

\bibitem{ASU-4}
 A. S. \"Ust\"unel: ``Variational calculation of Laplace transforms
 via entropy on Wiener space and applications''. 
J. Funct. Anal. {\bf{267}}, no. 8, 3058–3083, 2014.


\bibitem{W}
S. Watanabe:``The Japanese Contributions to Martingales''.
Journal électronique d’Histoire des Probabilités et de la Statistique/ Electronic Journal for 
History of Probability and Statistics . Vol.5, n°1. Juin/June 2009



\bibitem{Y-W} Y. Yamada, S. Watanabe:`` On the uniqueness of solutions
  of stochastic differential equations''. 
J. Math. Kyoto Univ. {\bf{11}} p. 155-167, 1971.




\end{thebibliography}
\end{document}